\documentclass[12pt]{article}

\usepackage{a4,amssymb,amsfonts,amsmath,graphicx,bbm}

\newcommand{\real}{\mathbbm{R}}
\newcommand{\complex}{\mathbbm{C}}
\newcommand{\nat}{\mathbbm{N}}
\newcommand{\ltwo}{\mathcal{L}^2}
\newcommand{\lonetime}{\mathcal{L}^1([0,\infty))}
\newcommand{\ltwotime}{\mathcal{L}^2([0,\infty))}
\newcommand{\ltwotimeupn}{\mathcal{L}^2([0,\infty))^n}
\newcommand{\ltwotimeneg}{\mathcal{L}^2((-\infty,0])}
\newcommand{\ltwopi}{\mathcal{L}^2(\Pi)}
\newcommand{\ltwoproduct}{\mathcal{L}^2(\Pi \times [0,\infty))}
\newcommand{\hinfty}{\mathcal{H}_{\infty}}

\newtheorem{theorem}{Theorem}
\newtheorem{lemma}{Lemma}
\newtheorem{definition}{Definition}

\begin{document}

\thispagestyle{plain}


\begin{center}
  \Large
  Sensitivity analysis of random linear dynamical systems
  using quadratic outputs

  \vspace{5mm}

  \vspace{0.5cm}   

{\large Roland~Pulch$\mbox{}^1$\footnote{corresponding author} and Akil Narayan$\mbox{}^2$}

\vspace{0.1cm}

\begin{small}

{$\mbox{}^1$Institute of Mathematics and Computer Science, 
University of Greifswald, \\
Walther-Rathenau-Str.~47, 
D-17489 Greifswald, Germany.}\\
Email: {\tt roland.pulch@uni-greifswald.de}

{$\mbox{}^2$Scientific Computing and Imaging Institute,
  Department of Mathematics, \\
  University of Utah, 
  72~Central Campus Dr., Salt Lake City, UT 84112,
  United States.}\\
  Email: {\tt akil@sci.utah.edu}

\end{small}
  
\end{center}

\bigskip\bigskip


\begin{center}
{Abstract}

\begin{tabular}{p{13cm}}
  In uncertainty quantification, a stochastic modelling is often applied,
  where parameters are substituted by random variables.
  We investigate linear dynamical systems of ordinary differential equations
  with a quantity of interest as output.
  Our objective is to analyse the sensitivity of the output
  with respect to the random variables.
  A variance-based approach generates partial variances and
  sensitivity indices.
  We expand the random output using the generalised polynomial chaos. 
  The stochastic Galerkin method yields a larger system of
  ordinary differential equations.
  The partial variances represent quadratic outputs of this system.
  We examine system norms of the stochastic Galerkin formulation
  to obtain sensitivity measures. 
  Furthermore, we apply model order reduction by balanced truncation,
  which allows for an efficient computation of the system norms with
  guaranteed error bounds.
  Numerical results are shown for a test example.

  \bigskip

\end{tabular}
\end{center}

\clearpage

\markright{R.~Pulch, A.~Narayan: Sensitivity analysis of random dynamical systems}


\section{Introduction}
Mathematical models of dynamical systems include physical or geometrical
parameters, which are often affected by uncertainties.
Thus the variability of the parameters has to be taken into account when using such models for predictive purposes.
In parametric uncertainty quantification (UQ), the varying parameters are often replaced
by random variables, see~\cite{sullivan,xiu-book}.
The solutions of a dynamical system therefore become random processes,
which encode the variability of the model output.

In this paper, we consider linear dynamical systems consisting of ordinary differential
equations (ODEs) that have random input parameters.
A time-dependent output is specified as a quantity of interest (QoI).
Stochastic modelling via UQ yields a random QoI, which must be examined.
Our main goal is to perform a sensitivity analysis of this random QoI with respect to the individual random parameters.
One of the most popular approaches for this strategy, Sobol indices,
is a technique for performing variance-based sensitivity analysis,
see~\cite{saltelli-etal,sobol,sobol-kucherenko}.
These sensitivity measures were employed to investigate dynamical systems
in different applications, for example, 
a motor model~\cite{haro-sandoval}, an epidemiological model~\cite{santonja}
and a (discretised) heat equation~\cite{soll-pulch}. 
Analysis via Sobol indices yields computationally tractable means of
determining total effect sensitivity indices,
which are defined by partial variances associated to each random parameter.
Previous related work investigates the sensitivity of the transfer function associated to random ODEs 
in the frequency domain~\cite{pulch-maten-augustin}.

Our approach for computing sensitivities utilises polynomial chaos (PC)
expansions;
the random QoI is expanded into a separation-of-variables series,
where each summand is the product of an orthogonal basis
polynomial with a time-dependent coefficient function, see~\cite{xiu-karniadakis}.
The stochastic Galerkin method with a PC ansatz yields a larger deterministic
linear dynamical system, whose outputs are approximations
to the time-dependent coefficients. 
In the previous work~\cite{pulch18}, the transfer function of the
stochastic Galerkin system was examined componentwise to
identify sparse approximations of the random QoI. We note that stochastic collocation approaches are likewise popular in parametric UQ settings, see~\cite{narayan_stochastic_2015}, but analysis for collocation schemes is more difficult, so we restrict our attention to Galerkin methods.

We investigate system norms of the stochastic Galerkin method
as sensitivity measures in the paper.
Each partial variance can be written as a quadratic output of the
stochastic Galerkin system. The first step in our approach transforms this system with a quadratic output into an 
equivalent linear dynamical system with multiple (linear)
outputs. 
The norm of the output of these linear dynamical systems represents a single
sensitivity coefficient for each random parameter.
In particular, we investigate the 
Hankel norm and a Hardy norm of the systems in our analysis.
The resulting sensitivity measures quantify the impact that
each random parameter has on the random QoI, and 
upper bounds on this variability are derived.
Thus the sensitivity analysis allows for the so-called
screening, see~\cite{ccs}, where random parameters with
sufficiently low impact are detected and remodelled as constants
to simplify the stochastic problem.

A major challenge in the above approach for sensitivity analysis is that the stochastic Galerkin dynamical system can be huge, and hence is computationally onerous to simulate.
Our contribution is application of model order reduction (MOR) strategies
to the high-dimensional stochastic Galerkin system,
which results in a much smaller reduced-order model that is much easier
to simulate. 
General information on MOR can be found in~\cite{antoulas,schilders};
in~\cite{pulch-maten}, the stochastic Galerkin method was reduced by
a Krylov subspace method.
In this paper, we use the technique of balanced truncation, where
a-priori error bounds are available for the MOR.
We provide analysis that bounds the proximity of the full-order model
sensitivity measures and the reduced-order model sensitivity measures,
so that our MOR approach is accompanied by a certifiable error estimate.


\section{Random linear dynamical systems}
We define the investigated problem in this section.

\subsection{Linear dynamical systems}
Let a time-invariant linear dynamical system be given in the form
\begin{equation} \label{ode}
  \begin{split}
    E(p) \dot{x}(t,p) &= A(p) x(t,p) + B(p) u(t) \\[1ex]
    y(t,p) &= C(p) x(t,p)
  \end{split}
\end{equation}
for $t \ge 0$
with inputs $u : [0,\infty) \rightarrow \real^{n_{\rm in}}$.
The matrices $A,E \in \real^{n \times n}$, 
$B \in \real^{n \times n_{\rm in}}$, $C \in \real^{n_{\rm out} \times n}$ 
depend on parameters $p \in \Pi \subseteq \real^{q}$.
We assume that the mass matrix~$E$ is non-singular for all $p \in \Pi$,
which implies a system of ODEs.
Thus both the state variables $x : [0,\infty) \times \Pi \rightarrow \real^n$ 
and the outputs 
$y : [0,\infty) \times \Pi \rightarrow \real^{n_{\rm out}}$ 
are also parameter-dependent. 
Without loss of generality, we restrict the problem to a single output
of the system ($n_{\rm out} = 1$).
Initial values
\begin{equation} \label{initial-values}
  x(0,p) = x_0(p)
\end{equation}
are predetermined by a function $x_0 : \Pi \rightarrow \real^n$.
The state variables solving the
initial value problem (\ref{ode}),(\ref{initial-values}) read as
\begin{equation} \label{ivp-solution}
  x(t,p) = {\rm e}^{E(p)^{-1}A(p) t} \left( x_0(p) +
  \int_0^t {\rm e}^{-E(p)^{-1}A(p) \tau} E(p)^{-1}B(p) u(\tau) \; {\rm d}\tau \right)
\end{equation}
for $t \ge 0$ and each $p \in \Pi$ including the matrix exponential.
However, this formula is not suitable for numerical evaluations.

Furthermore, we suppose that the system~(\ref{ode}) is always
asymptotically stable, i.e., all eigenvalues~$\lambda \in \complex$
satisfying the condition
$$ \det(\lambda E(p) - A(p)) = 0 $$
have a negative real part for each $p \in \Pi$.
Often the asymptotic stability of a model depends only on the topology
of the underlying configuration and not on the values of the
physical parameters, while still requiring positivity of some parameters.
In~\cite[p.~124]{riaza-tischendorf}, the authors show a sufficient condition
for linear electric circuits based on the topology.
In~\cite{hoagg-etal}, the asymptotic stability of a class of serial
mass-spring-damper systems is proven.

\subsection{Stochastic modelling}
If the parameters are affected by uncertainties,
then we use independent random variables on a probability space
to model the parameter variations.
The random variables are $p : \Omega \rightarrow \Pi$,
$\omega \mapsto p(\omega)$ with an event space~$\Omega$.
Let a joint probability density function $\rho : \Pi \rightarrow \real$
be given.
If a measurable function $f : \Pi \rightarrow \real$ depends on the
parameters, then its expected value reads as
\begin{equation} \label{expected-value}
  \mathbb{E} [f] = \int_{\Pi} f(p) \, \rho(p) \; {\rm d}p
\end{equation}
provided that the integral is finite.
The Hilbert space
\begin{equation} \label{l2-space}
  \ltwopi = \left\{ f : \Pi \rightarrow \real \; : \;
  f \; \mbox{measurable and} \; \mathbb{E} [f^2] < \infty \right\}
\end{equation}
is equipped with the inner product
\begin{equation} \label{innerproduct}
  \langle f , g \rangle =  \int_{\Pi} f(p) g(p) \, \rho(p) \; {\rm d}p
  \qquad \mbox{for} \;\; f,g \in \ltwopi .
\end{equation}
The induced norm reads as
$\| f \|_{\ltwopi} = \sqrt{\langle f,f \rangle}$.
We assume that the matrices of the linear dynamical system~(\ref{ode})
as well as the initial values~(\ref{initial-values})
depend continuously on the parameters.
The state variables are given by the formula~(\ref{ivp-solution}).
If the parameter domain~$\Pi$ is compact, then the state variables
$x(t,p)$ are uniformly bounded for all~$p$ and fixed time~$t$.
It follows that the state variables are located in
the Hilbert space~(\ref{l2-space}) pointwise for each~$t$.
If the parameter domain~$\Pi$ is infinite,
then integrability conditions have to be satisfied to guarantee
that the state variables belong to~(\ref{l2-space}) pointwise in time.
The output inherits the $\ltwopi$-property of the state variables if
the matrix~$C(p)$ is uniformly bounded for all~$p \in \Pi$.

Let a complete orthonormal system
$(\Phi_i)_{i \in \nat} \subset \ltwopi$
consisting of polynomials be given, i.e., 
$$ \langle \Phi_i , \Phi_j \rangle = \left\{
\begin{array}{ll}
  1 & \mbox{for} \; i=j , \\
  0 & \mbox{for} \; i \neq j . \\
\end{array} \right. $$
The theory of the generalised polynomial chaos (PC) implies a family of
orthonormal polynomials for each traditional random distribution,
see~\cite{xiu-book}.
The multivariate polynomials are the products of univariate
orthonormal polynomials.
The total degree of a multivariate polynomial is the sum of the
degrees in the univariate polynomials.
Let the polynomials be ordered such that
${\rm degree} (\Phi_i) \le {\rm degree} (\Phi_{i+1})$
is satisfied.
It follows that $\Phi_1$ is the unique constant polynomial.

Consequently, we expand the state variables as well as the output of the
system~(\ref{ode}) into
\begin{equation} \label{pc}
  x(t,p) = \sum_{i=1}^\infty v_i(t) \Phi_i(p)
  \qquad \mbox{and} \qquad
  y(t,p) = \sum_{i=1}^\infty w_i(t) \Phi_i(p)
\end{equation}
with coefficient functions $v_i : [0,\infty) \rightarrow \real^n$
and $w_i : [0,\infty) \rightarrow \real$.
The series converge in the norm of the Hilbert space~(\ref{l2-space}) and
pointwise in time.
The total variance of the random QoI becomes
\begin{equation} \label{total-variance}
  V(t) = \sum_{i=2}^\infty w_i(t)^2
\end{equation}
for each~$t \ge 0$.

\subsection{Stochastic Galerkin method}
In a numerical approximation, we include all polynomials up to a
total degree~$d$ described by the set of integers
\begin{equation} \label{indices-d}
  \mathcal{I}^d =
  \left \{ i \in \nat \; : \; {\rm degree}(\Phi_i) \le d \right\} .
\end{equation}
The cardinality is
$| \mathcal{I}^d |= \frac{(d+q)!}{d!q!}$, see~\cite[p.~65]{xiu-book}.
Hence the infinite summation in~(\ref{pc}) is restricted to
\begin{equation} \label{truncated-pc}
  x^{(d)}(t,p) = \sum_{i \in \mathcal{I}^d} v_i(t) \Phi_i(p)
  \qquad \mbox{and} \qquad
  y^{(d)}(t,p) = \sum_{i \in \mathcal{I}^d} w_i(t) \Phi_i(p) ,
\end{equation}
which represent truncated series.
Inserting the approximations~(\ref{truncated-pc}) into the
linear dynamical system~(\ref{ode}) yields a residual.
The Galerkin approach requires the residual to be orthogonal to
the space spanned by $\{ \Phi_i : i \in \mathcal{I}^d \}$
with respect to the inner product~(\ref{innerproduct}).
We obtain the larger linear system of ODEs
\begin{equation} \label{galerkin}
  \begin{split}
    \hat{E} \dot{\hat{v}}(t) &= \hat{A} \hat{v}(t) + \hat{B} u(t) \\[1ex]
    \hat{w} (t) &= \hat{C} \hat{v}(t)
  \end{split}
\end{equation}
with state variables~$\hat{v} : [0,\infty) \rightarrow \real^{mn}$
and the same inputs~$u$ as in~(\ref{ode}).
The system produces $m$~outputs~$\hat{w}$ with $m = | \mathcal{I}^d |$,
which represent approximations of the exact coefficient
functions in~(\ref{truncated-pc}).
Let
$$ S(p) = (\Phi_i(p)\Phi_j(p))_{i,j=1,\ldots,m} 
\quad \mbox{and} \quad
s(p) = (\Phi_i(p))_{i=1,\ldots,m}  $$
be auxiliary arrays.
The definition of the matrices
$\hat{A},\hat{E} \in \real^{mn \times mn}$, $\hat{B} \in \real^{mn \times n_{\rm in}}$,
$\hat{C} \in \real^{m \times mn}$
reads as
\begin{equation} \label{galerkin-matrices}
  \hat{A} = \mathbb{E} [S \otimes A] , \quad
  \hat{B} = \mathbb{E} [s \otimes B] , \quad
  \hat{C} = \mathbb{E} [S \otimes C] , \quad
  \hat{E} = \mathbb{E} [S \otimes E]
\end{equation}
using Kronecker products, where the
probabilistic integration~(\ref{expected-value}) is applied componentwise.
If the matrices $A,B,C,E$ consist of polynomials in the random variables,
then the matrices~(\ref{galerkin-matrices}) can be calculated analytically
for traditional probability distributions.
This property is an advantage in comparison to stochastic collocation
techniques, where a quadrature error or sampling error emerges. 
More details on the stochastic Galerkin method for linear dynamical systems
can be found in~\cite{pulch-maten}.

Initial values $\hat{v}(0) = \hat{v}_0$ have to be determined from
the initial values~(\ref{initial-values}) of the original system~(\ref{ode}).
If the initial values~(\ref{initial-values}) are identical to zero,
then the choice $\hat{v}(0) = 0$ is obvious.

The stochastic Galerkin system~(\ref{galerkin}) may be unstable
even if all original systems~(\ref{ode}) are asymptotically stable,
see~\cite{pulch-augustin}.
However, this loss of stability is rather seldom in practise.
Thus we assume that the linear dynamical system~(\ref{galerkin})
is asymptotically stable.


\section{Sensitivity analysis}
\label{sec:sensitivity}
We investigate the sensitivity of the random output from the
linear dynamical system~(\ref{ode}) with respect to the
individual random parameters.
Local sensitivity measures are based on partial derivatives with
respect to the parameters.
Alternatively, we are interested in global variance-based
sensitivity measures.

\subsection{Partial variances and sensitivity indices}
The Sobol indices provide a set of non-negative real numbers,
which describe the interaction of each subset of random variables,
see~\cite{sobol,sobol-kucherenko}.
Thus the number of Sobol indices is equal to the number of non-empty subsets,
i.e., $2^q-1$.
The Sobol indices yield the total effect sensitivity coefficients,
which represent a variance-based sensitivity measure.

We consider an equivalent definition of the total effect sensitivity
indices, which applies the PC expansion of a random-dependent
function, cf.~\cite{sudret}.
The random output of the system~(\ref{ode}) exhibits the
PC expansion~(\ref{pc}).
Let
\begin{equation} \label{indices-variable}
  \mathcal{I}_j = \left \{ i \in \nat \; : \;
  \Phi_i \; \mbox{is non-constant in} \; p_j \right\}
\end{equation}
be a set of integers for $j=1,\ldots,q$.
Considering all polynomials up to a total degree~$d$,
we obtain the intersection of~(\ref{indices-d}) and~(\ref{indices-variable})
\begin{equation} \label{indices-subset}
  \mathcal{I}_j^d = \mathcal{I}_j \cap \mathcal{I}^d
\end{equation}
for $j=1,\ldots,q$.
The partial variance associated to the $j$th random variable becomes
\begin{equation} \label{partial-var}
  V_j(t) = \sum_{i \in \mathcal{I}_j} w_i(t)^2
\end{equation}
for $j=1,\ldots,q$ and each $t \ge 0$.
It follows that $0 \le {V}_j(t) \le {V}(t)$ for all~$j$.
We assume that the total variance~(\ref{total-variance}) is positive
for all $t > 0$.
(It would be a rare exception if the variance reduces exactly to zero at
some positive time.
In this case, variations of the QoI vanish and a UQ becomes obsolete.)
The total effect sensitivity indices read as
\begin{equation} \label{total-sensitivity}
  S_j^{\rm T}(t) = \frac{V_j(t)}{V(t)} 
\end{equation}
for $j=1,\ldots,q$.
It holds that $0 \le S_j^{\rm T} \le 1$ for all~$j$ and
$1 \le S_1^{\rm T} + \cdots + S_q^{\rm T} \le q$ pointwise in time.
Often the sum of the sensitivity indices is close to one.

The stochastic Galerkin system~(\ref{galerkin}) yields approximations
of~(\ref{total-sensitivity}) by
\begin{equation} \label{sens-appr}
  \hat{S}_j^{\rm T}(t) = \frac{\hat{V}_j(t)}{\hat{V}(t)}
\end{equation}
with 
\begin{equation} \label{partial-var-appr}
  \hat{V}_j(t) = \sum_{i \in \mathcal{I}_j^d} \hat{w}_i(t)^2
\end{equation}
using the outputs $\hat{w}_1,\ldots,\hat{w}_m$
for the approximation of partial variances and, likewise,
for the approximation~$\hat{V}$ of the total variance~(\ref{total-variance}).

In our context, the sensitivity indices~(\ref{total-sensitivity})
are time-dependent functions,
which makes a sensitivity analysis more complicated.
The choice of the input signals and the initial values determine
these transient functions.
Alternatively, we want to obtain a single number as sensitivity measure
for each random variable, i.e., a finite set of $q$ non-negative
coefficients.

\subsection{Stochastic Galerkin system with quadratic outputs}
Let $\hat{w}^{(j)}$ be the vector-valued function consisting of the
components of $\hat{w}$ in $\mathcal{I}_j^d$ from~(\ref{indices-subset}).
It follows that
\begin{equation} \label{partial-var-2}
  \hat{V}_j(t) = \hat{w}^{(j)}(t)^\top \hat{w}^{(j)}(t) 
\end{equation}
for $j=1,\ldots,q$ due to~(\ref{partial-var-appr}).
We arrange an own output matrix for this partial variance now.
A square diagonal matrix is defined via
$$ D_j = {\rm diag} \left( d_1^{(j)} , \ldots , d_m^{(j)} \right)
\quad \mbox{with} \quad
  d_{\ell}^{(j)} = \left\{
  \begin{array}{rcl}
    1 & \mbox{if} \; \ell \in \mathcal{I}_j^d \\
    0 &  \mbox{if} \; \ell \notin \mathcal{I}_j^d \\
  \end{array} \right. $$
for each~$j=1,\ldots,q$.
The matrix-matrix product $\hat{C}_j' = D_j \hat{C}$ 
replaces rows of~$\hat{C}$ by zeros.
The rank of~$\hat{C}_j'$ is $k = | \mathcal{I}_j^d |$,
since the rank of~$\hat{C}$ is $m$.
We remove the rows identical to zero in $\hat{C}_j'$
to obtain a condensed matrix $\hat{C}_j \in \real^{k \times n}$.
Observing~(\ref{partial-var-2}), it follows that
\begin{equation} \label{variance-C}
  \hat{V}_j(t) = \hat{w}^{(j)}(t)^\top \hat{w}^{(j)}(t) =
  \hat{v}(t)^\top \hat{C}_j'^\top \hat{C}_j' \hat{v}(t) =
  \hat{v}(t)^\top \hat{C}_j^\top \hat{C}_j \hat{v}(t) .
\end{equation}
The partial variance of the $j$th random variable can be obtained
from a stochastic Galerkin system
\begin{equation} \label{galerkin-quadratic}
  \begin{split}
    \hat{E} \dot{\hat{v}}(t) &= \hat{A} \hat{v}(t) + \hat{B} u(t) \\[1ex]
    \hat{V}_j (t) &= \hat{v}(t)^\top \hat{C}_j^\top \hat{C}_j \hat{v}(t) ,
  \end{split}
\end{equation}
where a single quadratic output is defined by the symmetric positive
semi-definite matrix
\begin{equation} \label{matrices-quadratic}
  M_j = \hat{C}_j^\top \hat{C}_j
\end{equation}
for each $j=1,\ldots,q$.

Likewise, the total variance is given as the single quadratic output
$$ \hat{V} (t) = \hat{v}(t)^\top \hat{C}^\top \hat{C} \hat{v}(t) $$
from the differential equations in the
stochastic Galerkin system~(\ref{galerkin}).

In~\cite{beeumen-meerbergen,beeumen2012}, the authors propose an
approach for a general linear dynamical system with a quadratic output
specified by a symmetric positive semi-definite matrix.
A symmetric decomposition of this matrix is required like
the (pivoted) Cholesky factorisation, for example.
In our application, we already possess such a decomposition
via~(\ref{matrices-quadratic}).
The alternative stochastic Galerkin systems 
\begin{equation} \label{galerkin2}
  \begin{split}
    \hat{E} \dot{\hat{v}}(t) &= \hat{A} \hat{v}(t) + \hat{B} u(t) \\[1ex]
    \hat{z}_j (t) &= \hat{C}_j \hat{v}(t) 
  \end{split}
\end{equation}
are arranged for $j=1,\ldots,q$
with multiple outputs ($n_{\rm out} = k = | \mathcal{I}_j^d |$).
It holds that $\hat{V}_j = \hat{z}_j^\top \hat{z}_j$ for the
quadratic output of~(\ref{galerkin-quadratic})
due to~(\ref{variance-C}).

\subsection{System norms}
\label{sec:system-norms}
Given a vector-valued function $x \in \ltwotimeupn$
with $n$~components, its signal norm reads as
$$ \| x \|_{\ltwotime} = \sqrt{ \int_0^\infty x(t)^\top x(t) \; {\rm d}t } \, . $$
Several norms are defined for a general linear dynamical system.
Let $\hat{\Sigma}$ and $\hat{\Sigma}_1,\ldots,\hat{\Sigma}_q$ represent
the stochastic Galerkin systems~(\ref{galerkin}) and~(\ref{galerkin2}),
respectively.
The Hankel norm denotes a bound on the signal norm of future outputs
generated by inputs in the past (input is zero for $t \ge 0$).
This norm exhibits the formula, cf.~\cite[p.~135]{antoulas},
\begin{equation} \label{hankel-norm}
  \| \hat{\Sigma} \|_{\rm H} =
  \sup
  \left\{ \frac{\| \hat{w} \|_{\ltwotime}}{\| u \|_{\ltwotimeneg}} \; : \;
  u \in \ltwotimeneg^{n_{\rm in}} , \; u \neq 0 \right\} .
\end{equation}
In addition, the Hankel norm is equal to the maximum Hankel singular
value of the system, which will be specified in
Section~\ref{sec:balanced}.
The system norm directly describes the mapping from inputs to outputs
provided that the initial values are zero, i.e., 
\begin{equation} \label{l2-norm-system}
  \| \hat{\Sigma} \|_{\ltwo} =
  \sup
  \left\{ \frac{\| \hat{w} \|_{\ltwotime}}{\| u \|_{\ltwotime}} \; : \;
  u \in \ltwotime^{n_{\rm in}} , \; u \neq 0 , \; \hat{v}(0) = 0 \right\} .
\end{equation}
Let $\hat{G} : \complex \backslash \Xi \rightarrow \complex^{m \times n_{\rm in}}$
be the transfer function of the stochastic Galerkin system~(\ref{galerkin})
in the frequency domain.
Therein,  $\Xi \subset \complex$ is a finite set of poles.
The magnitude of a transfer function can be characterised by the
Hardy norms $\| \cdot \|_{\mathcal{H}_2}$ and $\| \cdot \|_{\hinfty}$.
The $\hinfty$-norm of a transfer function coincides with its
$\mathcal{L}^{\infty}$-norm, which is a frequency domain Lebesgue norm,
provided that the system is asymptotically stable.
It follows that, see~\cite[p.~149]{antoulas},
\begin{equation} \label{hinf-norm}
  \| \hat{G} \|_{\hinfty} = \sup_{\omega \in \real}
  \; \sigma_{\max} ( \hat{G}({\rm i}\omega) ) ,
\end{equation}
where $\sigma_{\max}$ is the maximum singular value of the
matrix-valued transfer function and ${\rm i} = \sqrt{-1}$.
The norms satisfy the relations
\begin{equation} \label{norm-relation}
  \| \hat{\Sigma} \|_{\rm H} \le \| \hat{\Sigma} \|_{\ltwo}
  = \| \hat{G} \|_{\hinfty} ,
\end{equation}
which are valid for any asymptotically stable system of ODEs.

\subsection{Sensitivity measures}
We use the norms of the stochastic Galerkin systems~(\ref{galerkin2})
from Section~\ref{sec:system-norms} to indicate the sensitivity of
the random output in the original system~(\ref{ode}).
In contrast, the $\mathcal{H}_2$-norms and
$\hinfty$-norms of the separate components of the
transfer function~$\hat{G}$ were analysed to identify a
sparse representation of the random output in~\cite{pulch18}.

The following two definitions introduce the key data
in our investigations.

\begin{definition} \label{def:measure-hankel}
  The sensitivity coefficients $\hat{\eta}_j$ are the real values
  \begin{equation} \label{sens-measure-hankel}
    \hat{\eta}_j = \| \hat{\Sigma}_j \|_{\rm H}
  \end{equation}
  for $j=1,\ldots,q$ using the Hankel norms~(\ref{hankel-norm}) of the
  stochastic Galerkin systems~(\ref{galerkin2}).
\end{definition}

\begin{definition} \label{def:measure}
  The sensitivity coefficients $\hat{\theta}_j$ are the real values
  \begin{equation} \label{sens-measure}
    \hat{\theta}_j = \| \hat{\Sigma}_j \|_{\ltwo} = \| \hat{G}_j \|_{\hinfty}
  \end{equation}
  for $j=1,\ldots,q$ using the system norms (\ref{l2-norm-system}) of the
  stochastic Galerkin systems~(\ref{galerkin2}).
\end{definition}

The general relation~(\ref{norm-relation}) shows that
\begin{equation} \label{eta-le-theta}
  \hat{\eta}_j \le \hat{\theta}_j \qquad \mbox{for each} \;\; j=1,\ldots,q .
\end{equation}
The sensitivity measure~(\ref{sens-measure}) is more important
than~(\ref{sens-measure-hankel}), because~(\ref{sens-measure-hankel})
does not directly relate to the input-output mapping for $t \ge 0$.
However, the computation of an $\hinfty$-norm is more costly than
the calculation of a Hankel norm.
We obtain an inequality between the system norms of the Galerkin
system~(\ref{galerkin}) and the modified Galerkin systems~(\ref{galerkin2}).

\begin{lemma} \label{lemma:partial-norm}
  Let $\hat{\Sigma}$ and $\hat{\Sigma}_1,\ldots,\hat{\Sigma}_q$ 
  represent the linear dynamical systems~(\ref{galerkin})
  and~(\ref{galerkin2}), respectively.
  It holds that $\| \hat{\Sigma}_j \|_{\ltwo} \le \| \hat{\Sigma} \|_{\ltwo}$
  for $j=1,\ldots,q$
  in the system norm~(\ref{l2-norm-system}).
\end{lemma}

Proof: \nopagebreak

The outputs of the system~(\ref{galerkin}) and the systems~(\ref{galerkin2})
are~$\hat{w}$ and $\hat{z}_1,\ldots,\hat{z}_q$, respectively.
We obtain the estimate
\begin{align*}
  \| \hat{w} \|_{\ltwotime}^2 &=
\int_0^\infty \sum_{i \in \mathcal{I}^d} \hat{w}_i(t)^2 \; {\rm d}t
\ge \int_0^\infty \sum_{i \in \mathcal{I}_j^d} \hat{w}_i(t)^2 \; {\rm d}t
=  \int_0^\infty \hat{z}_j(t)^\top \hat{z}_j(t) \; {\rm d}t \\
& =  \| \hat{z}_j \|_{\ltwotime}^2
\end{align*}
for $j=1,\ldots,q$.
It follows that $\| \hat{z}_j \|_{\ltwotime} \le \| \hat{w} \|_{\ltwotime}$
if the same input is supplied to all systems.
Observing~(\ref{l2-norm-system}), we arrange the supremum and
achieve the statement.
\hfill $\Box$

\medskip

We show a bound on the norm of partial variances in the time domain,
which yields an interpretation of the sensitivity measure
in Definition~\ref{def:measure}.

\begin{lemma} \label{lemma:partial-norm-bound}
The stochastic Galerkin systems~(\ref{galerkin2}) with initial values zero
  yield approximations of the partial variances,
  which satisfy the bounds
  $$ \| \hat{V}_j \|_{\lonetime} \le
      \hat{\theta}_j^2 \, \| u \|_{\ltwotime}^2 $$
  for $j=1,\ldots,q$ including the
  sensitivity coefficients~(\ref{sens-measure}).
\end{lemma}

Proof:

We obtain
\begin{align*}
  \| \hat{V}_j \|_{\mathcal{L}^1[0,\infty)} & =
  \int_0^\infty  \hat{V}_j(t) \; {\rm d}t =
  \int_0^\infty \hat{z}_j(t)^\top \hat{z}_j(t) \; {\rm d}t
  =  \| \hat{z}_j \|_{\ltwotime}^2 \\[1ex]
  & \le \| \hat{\Sigma}_j \|_{\ltwo}^2 \| u \|_{\ltwotime}^2 =
  \hat{\theta}_j^2 \| u \|_{\ltwotime}^2
\end{align*}
due to~(\ref{l2-norm-system}) and~(\ref{sens-measure}).
\hfill $\Box$

\medskip

The definition of the total effect sensitivity
indices~(\ref{total-sensitivity}) motivates to establish 
relative sensitivity measures
$$
\hat{\eta}_j^{\rm rel} =
\frac{\| \hat{\Sigma}_j \|_{\rm H}}{\| \hat{\Sigma} \|_{\rm H}}
\qquad \mbox{as well as} \qquad
\hat{\theta}_j^{\rm rel} =
\frac{\| \hat{\Sigma}_j \|_{\ltwo}}{\| \hat{\Sigma} \|_{\ltwo}} $$ 
for $j=1,\ldots,q$.
Lemma~\ref{lemma:partial-norm} guarantees
$\hat{\theta}_j^{\rm rel} \le 1$, whereas $\hat{\eta}_j^{\rm rel}$ may be larger than one.
However, since the denominators are identical for $j=1,\ldots,q$,
we can investigate the normalised quantities
\begin{align} 
  \hat{\eta}_j^* &=
  \frac{\hat{\eta}_j}{\hat{\eta}_1 + \cdots + \hat{\eta}_q} ,
  \label{sensitivity-normalised1} \\
  \hat{\theta}_j^* &=
  \frac{\hat{\theta}_j}{\hat{\theta}_1 + \cdots + \hat{\theta}_q}
  \label{sensitivity-normalised2} 
\end{align}
as well.
Now it holds that $0 \le \hat{\eta}_j^* , \hat{\theta}_j^* \le 1$
for all $j=1,\ldots,q$.

\subsection{Error bound in screening}
The sensitivity measures identify the random variables
with high impact and low impact.
Random variables with low impact can be remodelled by constants
to decrease the dimensionality of the random parameter space.
This strategy is called screening, see~\cite{ccs}, or
freezing of unessential random variables, see~\cite{sobol}.
We derive a bound on the error of this screening.

\begin{theorem} \label{thm:sensitivity}
  Let the random variables be partitioned into $p=(a,b)$ with
  $a \in \Pi_a \subseteq \real^{q'}$, $b \in \Pi_b \subseteq \real^{q-q'}$.
  Let $\bar{y}(t,p) = y(t,a,\bar{b})$ with a constant~$\bar{b} \in \Pi_b$.
  For each $\delta \in (0,1)$, there is a set
  $B_{\delta} \subset \Pi_b$ of probability larger $1-\delta$
  such that
  $$ \left\| y - \bar{y} \right\|_{\ltwoproduct}^2 <
  \left( 1 + {\textstyle \frac{1}{\delta}} \right)
  \sum_{j=q'+1}^q \| V_j \|_{\lonetime} $$
  for all $\bar{b} \in B_{\delta}$
  in the norm of the product space $\Pi \times [0,\infty)$.
\end{theorem}

Proof: \nopagebreak

The analysis in~\cite{sobol} shows that
$$  \int_{\Pi} \left( y(t,a,b) - y(t,a,\bar{b}) \right)^2 \rho(p) \; {\rm d}p
< \left( 1 +  {\textstyle \frac{1}{\delta}} \right) \sum_{j=q'+1}^q V_j(t) $$
pointwise for each~$t \ge 0$.
Time integration in $[0,\infty)$ yields
  \begin{align*}
    \left\| y - \bar{y} \right\|_{\ltwoproduct}^2 &=
    \int_0^\infty \int_{\Pi} \left( y(t,p)-\bar{y}(t,p) \right)^2 \rho(p)
        \; {\rm d}p \; {\rm d}t \\
    &< \left( 1 +  {\textstyle \frac{1}{\delta}} \right)
    \int_0^\infty \hspace{-3mm} \sum_{j=q'+1}^q V_j(t) \; {\rm d}t 
    = \left( 1 +  {\textstyle \frac{1}{\delta}} \right)
   \sum_{j=q'+1}^q \| V_j \|_{\lonetime} ,
  \end{align*}
  because the partial variances are non-negative.
   \hfill $\Box$

   \medskip

The stochastic Galerkin systems generate approximations of the
partial variances~(\ref{partial-var}).
Assuming that $\hat{V}_j - V_j$ is sufficiently small
for $j=q'+1,\ldots,q$,
Lemma~\ref{lemma:partial-norm-bound} and Theorem~\ref{thm:sensitivity}
imply the approximate bound
$$  \left\| y - \bar{y} \right\|_{\ltwoproduct} <
\sqrt{ {\textstyle \left( 1 + \frac{1}{\delta} \right)}
  \sum_{j=q'+1}^q \hat{\theta}_j^2 \;\mbox{}} \;
\| u \|_{\ltwotime} . $$
This inequality shows that the random variables~$b$ can be replaced
by constants if the associated sensitivity
measures~(\ref{sens-measure}) are sufficiently small.


\section{Model order reduction}
We examine the potential for a cheap computation of the
sensitivity measures in Definition~\ref{def:measure-hankel}
and Definition~\ref{def:measure} using methods of MOR.

\subsection{Balanced truncation}
\label{sec:balanced}
We apply an MOR to the stochastic Galerkin system~(\ref{galerkin})
using the technique of balanced truncation, see~\cite{antoulas}.
This MOR method is applicable to all asymptotically stable linear systems
of ODEs.
The controllability Gramian~$P$ and the observability Gramian~$Q$
are given by the solutions of the Lyapunov equations
\begin{align} 
  \hat{A} P \hat{E}^\top + \hat{E} P \hat{A}^\top + \hat{B} \hat{B}^\top &= 0 ,
  \label{control-gramian} \\[0.5ex]
  \hat{A}^\top Q \hat{E} + \hat{E}^\top Q \hat{A} + \hat{C}^\top \hat{C} &= 0 .
  \label{observ-gramian}
\end{align}
We require symmetric decompositions
$P = L_P {L_P}^\top$ and $Q = L_Q {L_Q}^\top$.
For example, the Cholesky factorisation can be directly computed
without solving for $P$ or $Q$ first, see~\cite{hammarling}.
The singular value decomposition (SVD) 
\begin{equation} \label{svd}
  U_{\rm l} S {U_{\rm r}}^\top = {L_Q}^\top \hat{E} L_P
\end{equation}
is calculated with a diagonal matrix
$S = {\rm diag}(\sigma_1, \ldots , \sigma_{mn})$ and
orthogonal matrices $U_{\rm l},U_{\rm r} \in \real^{mn \times mn}$.
The diagonal matrix includes the Hankel singular values
in descending order ($\sigma_{\ell} \ge \sigma_{\ell+1}$).

Any reduced dimension $r < mn$ can be chosen.
Let $\tilde{S} \in \real^{r \times r}$ be the diagonal matrix including
the $r$ largest singular values.
Let $\tilde{U}_{\rm l},\tilde{U}_{\rm r} \in \real^{mn \times r}$ be the
matrices consisting of first $r$ columns in $U_{\rm l},U_{\rm r}$. 
We obtain projection matrices $T_{\rm l},T_{\rm r} \in \real^{mn \times r}$ of
full (column) rank by
$$ T_{\rm l} =  L_Q \tilde{U}_{\rm l} \tilde{S}^{-\frac{1}{2}} 
\qquad \mbox{and} \qquad
T_{\rm r} = L_P \tilde{U}_{\rm r} \tilde{S}^{-\frac{1}{2}} . $$
Alternatively, just the $r$ dominant singular values and their
singular vectors have to be computed to obtain these projection matrices.

The reduced matrices read as
\begin{equation} \label{reduced-matrices}
  \tilde{A} = T_{\rm l}^\top \hat{A} T_{\rm r} , \quad
  \tilde{B} = T_{\rm l}^\top \hat{B} , \quad
  \tilde{C} = \hat{C} T_{\rm r} , \quad
  \tilde{E} = T_{\rm l}^\top \hat{E} T_{\rm r} .
\end{equation}
The matrices produce a linear dynamical system of the same form but
smaller dimension~$r$.
The reduced mass matrix~$\tilde{E}$ is non-singular.
Let $\tilde{w}$ be its multiple output,
where the number of outputs is still~$m$.

We use the same projection matrices for reductions of the
linear dynamical systems~(\ref{galerkin2}).
Just the matrices
$$ \tilde{C}_j = \hat{C}_j T_{\rm r} $$
are changed for $j=1,\ldots,q$.
Moreover, the balanced truncation technique guarantees that each
reduced-order model (ROM) inherits the asymptotic stability of the
full-order model (FOM).

Direct methods of linear algebra yield the numerical solution
of the Lyapunov equations (\ref{control-gramian}),(\ref{observ-gramian}),
see~\cite{hammarling}.
Alternatively, approximate methods or iteration schemes are available
like Krylov subspace techniques and the alternating direction implicit (ADI)
method, see~\cite{druskin}.
Such methods produce low-rank approximations of the Cholesky factors. 
However, the efficiency is often critical
if the matrices $\hat{B}$ or $\hat{C}$ exhibit a high rank.
The number of inputs typically satisfies $n_{\rm in} < n \ll mn$
and~(\ref{control-gramian}) can be solved efficiently.
In contrast, it holds that ${\rm rank}(\hat{C})=m$
and thus the rank is relatively large for moderate~$n$,
where the approximate solution of~(\ref{observ-gramian})
becomes critical.

\subsection{Error bound for sensitivity measures}
We show an a-priori error bound on our sensitivity measures from
Definition~\ref{def:measure-hankel} and Definition~\ref{def:measure}.

\begin{theorem} \label{thm:error}
  Let $\hat{\Sigma}_j$ be the linear stochastic Galerkin
  system~(\ref{galerkin2}) and $\tilde{\Sigma}_j$ be its ROM of dimension~$r$
  obtained by balanced truncation of the complete
  stochastic Galerkin system~(\ref{galerkin}).
  The error of the sensitivity measures~(\ref{sens-measure-hankel})
  and~(\ref{sens-measure}) satisfy the bound
  \begin{equation} \label{error-bound}
    \max \left\{
  \big| \hat{\eta}_j - \tilde{\eta}_j \big| ,
  \big| \hat{\theta}_j - \tilde{\theta}_j \big|
  \right\} 
  \le 2 \sum_{k=r+1}^{mn} \sigma_k
  \end{equation}
  uniformly for all $j=1,\ldots,q$.
  The Hankel singular values
  $\sigma_1 \ge \sigma_2 \ge \cdots \ge \sigma_{mn}$
  belong to the stochastic Galerkin system~(\ref{galerkin}).
\end{theorem}

Proof:

It holds that $\hat{\eta}_j = \| \hat{\Sigma}_j \|_{\rm H}$,
$\hat{\theta}_j = \| \hat{\Sigma}_j \|_{\ltwo}$ and
$\tilde{\eta}_j = \| \tilde{\Sigma}_j \|_{\rm H}$,
$\tilde{\theta}_j = \| \tilde{\Sigma}_j \|_{\ltwo}$
for the sensitivity measures from FOM and ROM, respectively.
The reverse triangle inequality of the norms yields
\begin{align*}
  \big| \hat{\eta}_j - \tilde{\eta}_j \big| & =
  \left| \| \hat{\Sigma}_j \|_{\rm H} - \| \tilde{\Sigma}_j \|_{\rm H} \right|
  \le \left\| \hat{\Sigma}_j - \tilde{\Sigma}_j \right\|_{\rm H} , \\
  \big| \hat{\theta}_j - \tilde{\theta}_j \big| & =
  \left| \| \hat{\Sigma}_j \|_{\ltwo} - \| \tilde{\Sigma}_j \|_{\ltwo} \right|
  \le \left\| \hat{\Sigma}_j - \tilde{\Sigma}_j \right\|_{\ltwo} .
\end{align*}
Let $\hat{G}$ and $\tilde{G}$ denote the transfer functions of
FOM and ROM, respectively.
Lemma~\ref{lemma:partial-norm} together with the general
relation~(\ref{norm-relation}) show that
$$ \left\| \hat{\Sigma}_j - \tilde{\Sigma}_j \right\|_{\rm H} \le
\left\| \hat{\Sigma}_j - \tilde{\Sigma}_j \right\|_{\ltwo} \le
\left\| \hat{\Sigma} - \tilde{\Sigma} \right\|_{\ltwo} =
\left\| \hat{G} - \tilde{G} \right\|_{\mathcal{H}_\infty} $$
uniformly for~$j=1,\ldots,q$, because the involved systems are
asymptotically stable.
Using balanced truncation, the error estimate
$$ \left\| \hat{G} - \tilde{G} \right\|_{\mathcal{H}_\infty} \le
2 \sum_{k=r+1}^{mn} \sigma_k $$
holds true in the $\hinfty$-norm~(\ref{hinf-norm}),
see~\cite[p.~212]{antoulas}.
\hfill $\Box$

\medskip

Theorem~\ref{thm:error} demonstrates that the approximation of
the sensitivity measures inherits the error bound of the
balanced truncation MOR.


\section{Illustrative example}
We apply the sensitivity analysis from Section~\ref{sec:sensitivity}
to a test example now.
We computed on a FUJITSU Esprimo P920
Intel(R) Core(TM) i5-4570 CPU with 3.20 GHz (4~cores)
and operation system Microsoft Windows~7.
The software package MATLAB~\cite{matlab2018} (version R2018a)
was used for all computations.

\subsection{Mass-spring-damper system}
A general linear mass-spring-damper configuration is described by a
se\-cond-order system of ODEs
\begin{equation} \label{second-order-ode}
  M \ddot{z}(t) + D \dot{z}(t) + K z(t) = f(t)
\end{equation}
for the positions $z : [0,\infty) \rightarrow \real^{\ell}$
with mass matrix~$M$, damping matrix~$D$, stiffness matrix~$K$
and an excitation~$f: [0,\infty) \rightarrow \real^{\ell}$.

\begin{figure}
\centering
\includegraphics[width=6cm]{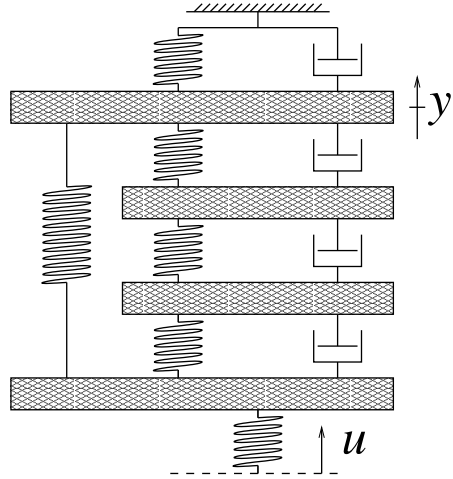}
\caption{Mass-spring-damper configuration.}
\label{fig:msd}       
\end{figure}

We examine a mass-spring-damper system introduced in~\cite{lohmann-eid},
which is depicted in Figure~\ref{fig:msd}.
This mechanical configuration incorporates $q=14$ positive
physical parameters: 
4~masses, 6~spring constants and 4~damping constants.
The single input is an excitation at the bottom spring,
whereas the single output is the position of the top mass.
The mathematical modelling yields the symmetric matrices
$M = {\rm diag}(m_1,m_2,m_3,m_4)$,
$$ D = \begin{pmatrix}
  d_1 & -d_1 & 0 & 0 \\
  -d_1 & d_1+d_2 & -d_2 & 0 \\
  0 & -d_2 & d_2+d_3 & -d_3 \\
  0 & 0 & -d_3 & d_3+d_4 \\
\end{pmatrix} , $$
$$ K = \begin{pmatrix}
  k_1+k_2+k_5 & -k_2 & 0 & -k_5 \\
  -k_2 & k_2+k_3 & -k_3 & 0 \\
  0 & -k_3 & k_3+k_4 & -k_4 \\
  -k_5 & 0 & -k_4 & k_4+k_5+k_6 \\
\end{pmatrix} $$
together with the excitation $f(t) = (k_1 u(t),0,0,0)^\top$.
Obviously, if all parameters are positive,
then $M$~is positive definite and $D$ as well as $K$ are
positive semi-definite.
A symbolic Cholesky decomposition is possible for both $D$ and~$K$. 
Thus these matrices are also positive definite.
It follows that the linear dynamical system is always
asymptotically stable, see~\cite[p.~103]{inman}.

\begin{figure}
\centering
\includegraphics[width=6.5cm]{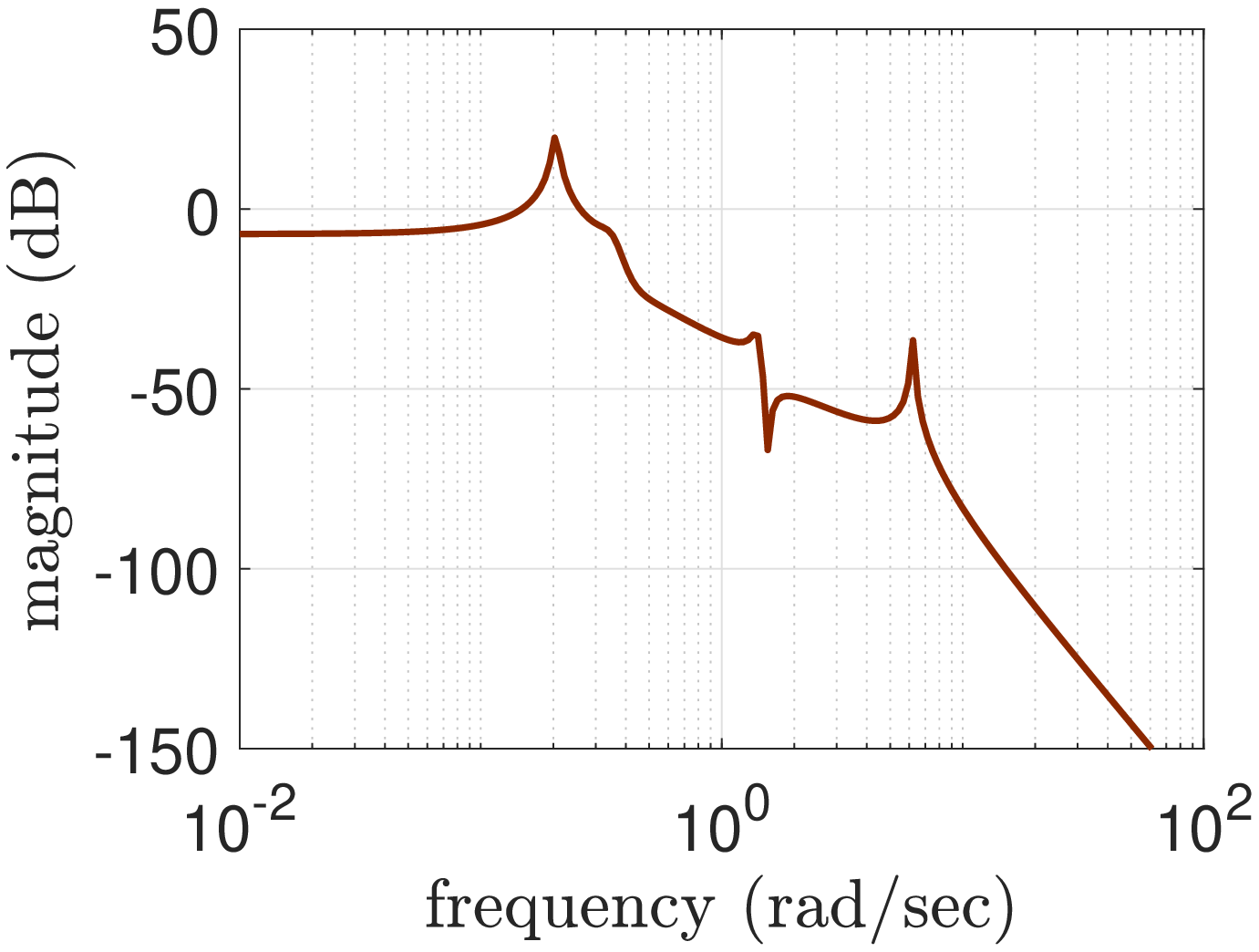}
\hspace{5mm}
\includegraphics[width=6.5cm]{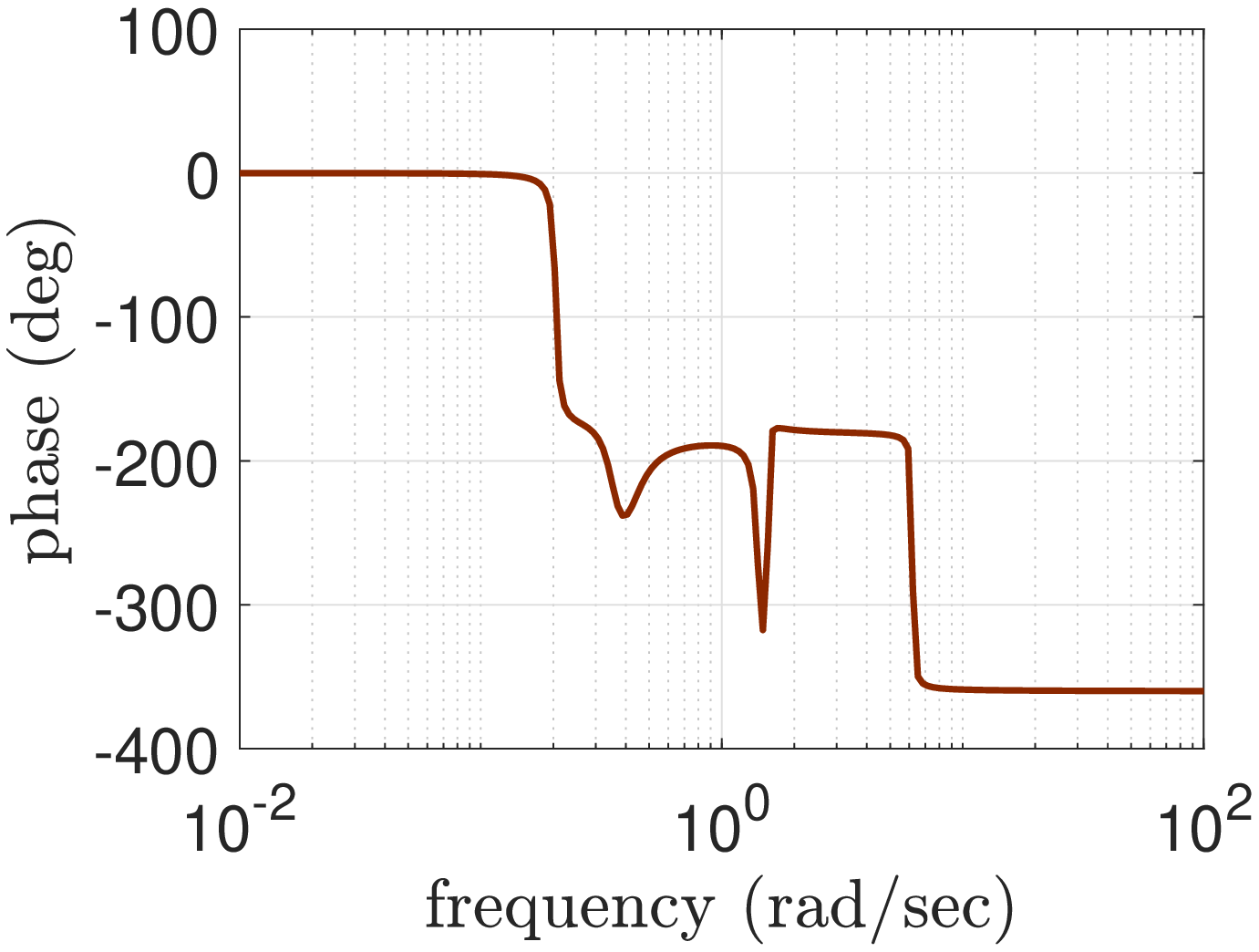}
\caption{Bode plot of mass-spring-damper system for
  deterministic parameters.}
\label{fig:bode}       
\end{figure}

The system~(\ref{second-order-ode}) is transformed into an equivalent
first-order system~(\ref{ode}) with dimension $n=8$.
Figure~\ref{fig:bode} shows the Bode plot of the linear dynamical system
in the case of a constant choice of the parameters from~\cite{lohmann-eid},
i.e., $m_1=1$, $m_2=5$, $m_3=25$, $m_4=125$,
$k_1 = 27$, $k_2 = 9$, $k_3 = 3$, $k_4 = 1$, $k_5 = 2$, $k_6 = 3$,
$d_1 = 0.1$, $d_2 = 0.4$, $d_3 = 1.6$, $d_4 = 1$.
This test example was also investigated in~\cite{pulch18,pulch-naco}.

\subsection{Stochastic modelling and Galerkin method}
We replace the parameters by independent random variables with
uniform probability distributions, which vary 10\% around the
used constant choice of the parameters.
Consequently, the PC expansions~(\ref{pc}) include
the (multivariate) Legendre polynomials.
We arrange finite sums~(\ref{truncated-pc}) with all polynomials
up to total degree~$d=3$.
The stochastic Galerkin system~(\ref{galerkin}) exhibits
the state space dimension $mn=5440$ and
the number of outputs becomes $m=680$.
The matrices of the Galerkin system can be computed analytically,
because the entries of matrices in the underlying
system~(\ref{second-order-ode})
are linear functions of the physical parameters.
Hence we calculate these matrices exactly except for round-off errors.
Furthermore, the stochastic Galerkin system is asymptotically stable,
because the maximum real part of the associated eigenvalues
is computed to $-0.0048$ and thus negative.

\subsection{Numerical results for sensitivity measures}
We computed all Hankel norms by direct methods of linear algebra.
In contrast, the $\hinfty$-norm could only be calculated by approximation.
We did not determine the $\hinfty$-norm of the system~(\ref{galerkin}),
because the computational effort would be too large
due to the high dimensionality.

\begin{figure}
\centering
\includegraphics[width=6.5cm]{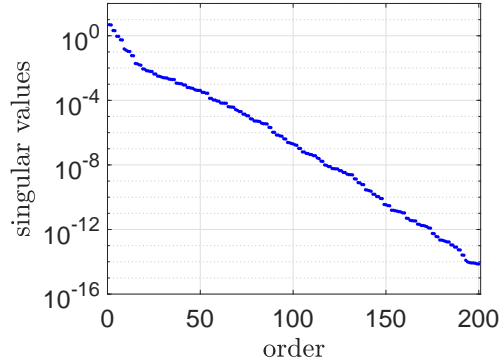}
\caption{Hankel singular values of the stochastic Galerkin
  system~(\ref{galerkin}) for the mass-spring-damper example.}
\label{fig:singularvalues}       
\end{figure}

We generate ROMs of the FOM~(\ref{galerkin}) by the
balanced truncation technique from Section~\ref{sec:balanced}.
A direct method of linear algebra yields the Cholesky factors of
the solutions satisfying the
Lyapunov equations~(\ref{control-gramian}),(\ref{observ-gramian}).
More details on the numerical method are given in the following subsection. 
Figure~\ref{fig:singularvalues} illustrates the dominant
Hankel singular values in the decomposition~(\ref{svd}).
The largest Hankel singular value coincides with the
Hankel norm $\| \hat{\Sigma} \|_{\rm H} = 4.84$.

\begin{figure}
\centering
\includegraphics[width=6.5cm]{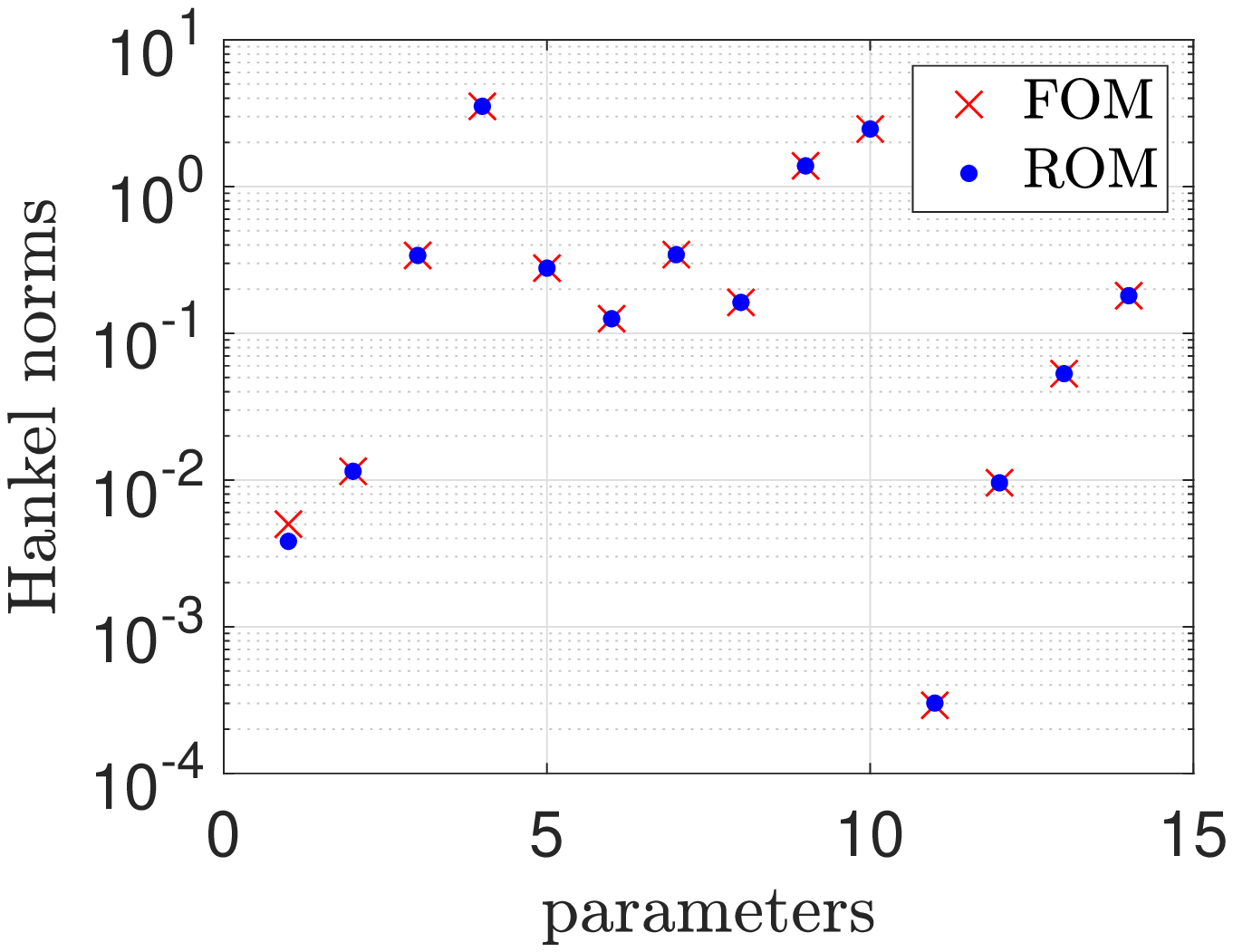}
\hspace{5mm}
\includegraphics[width=6.5cm]{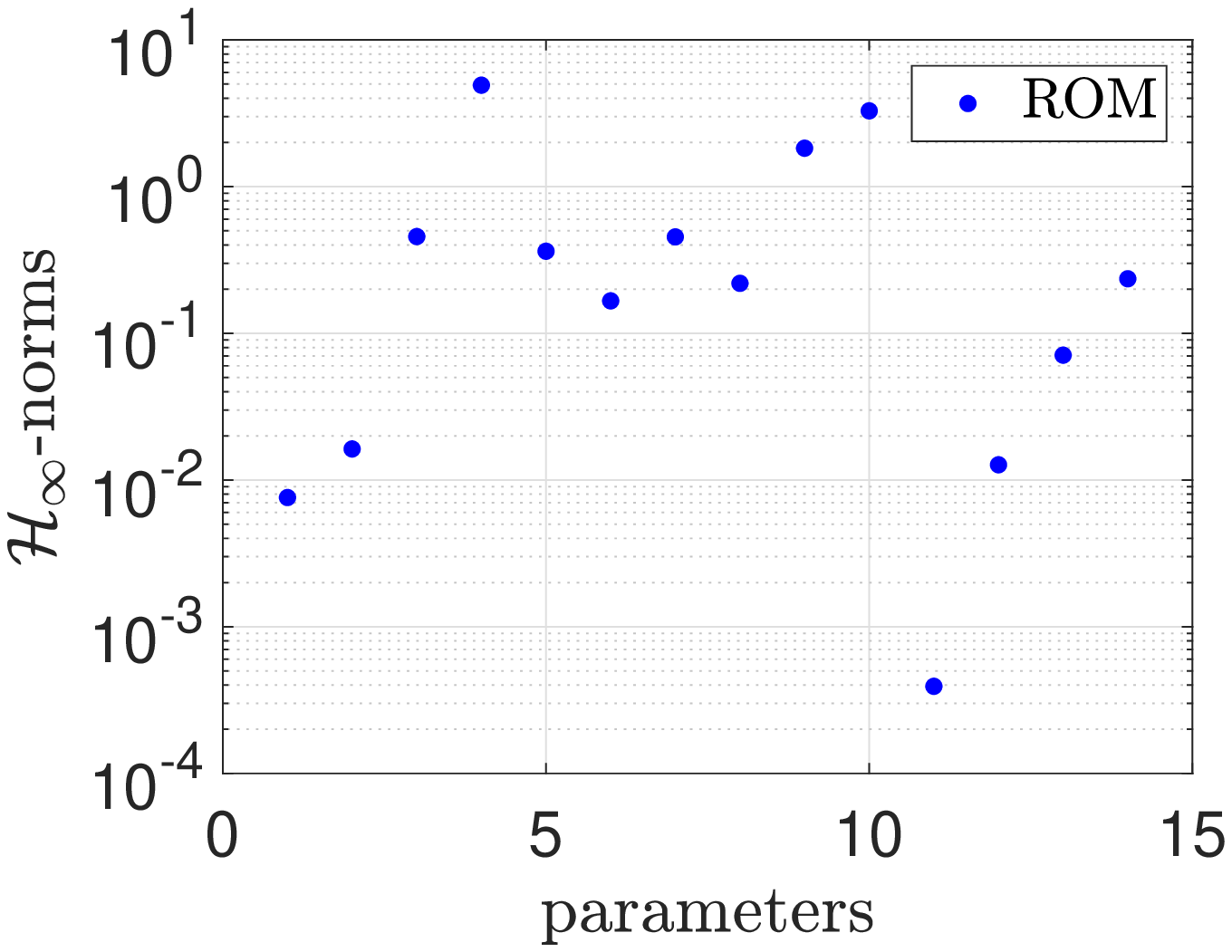}
\caption{Sensitivity coefficients~(\ref{sens-measure-hankel})
  from Hankel norms (left) and sensitivity coefficients~(\ref{sens-measure})
  from Hardy norms (right) with respect to the random parameters
  (1-4: masses, 5-10: springs, 11-14: dampers).}
\label{fig:sensitivities}       
\end{figure}

Now we analyse the sensitivities.
On the one hand, the FOMs~(\ref{galerkin2}) yield the sensitivity
measures~(\ref{sens-measure-hankel}).
On the other hand, the ROMs of dimension $r=25$ produce approximations
of the sensitivity coefficients~(\ref{sens-measure-hankel})
as well as~(\ref{sens-measure}). 
Figure~\ref{fig:sensitivities} depicts the results,
which demonstrate a good agreement between the sensitivity measures from
FOM and ROM.
Concerning the inequality~(\ref{eta-le-theta}),
the $\hinfty$-norms are just slightly larger than the Hankel norms.
Furthermore, the system norm of the ROM for the stochastic Galerkin
method~(\ref{galerkin}) is
$\| \tilde{\Sigma} \|_{\ltwo} = \| \tilde{G} \|_{\hinfty} = 7.75$.

\begin{figure}
\centering
\includegraphics[width=5cm]{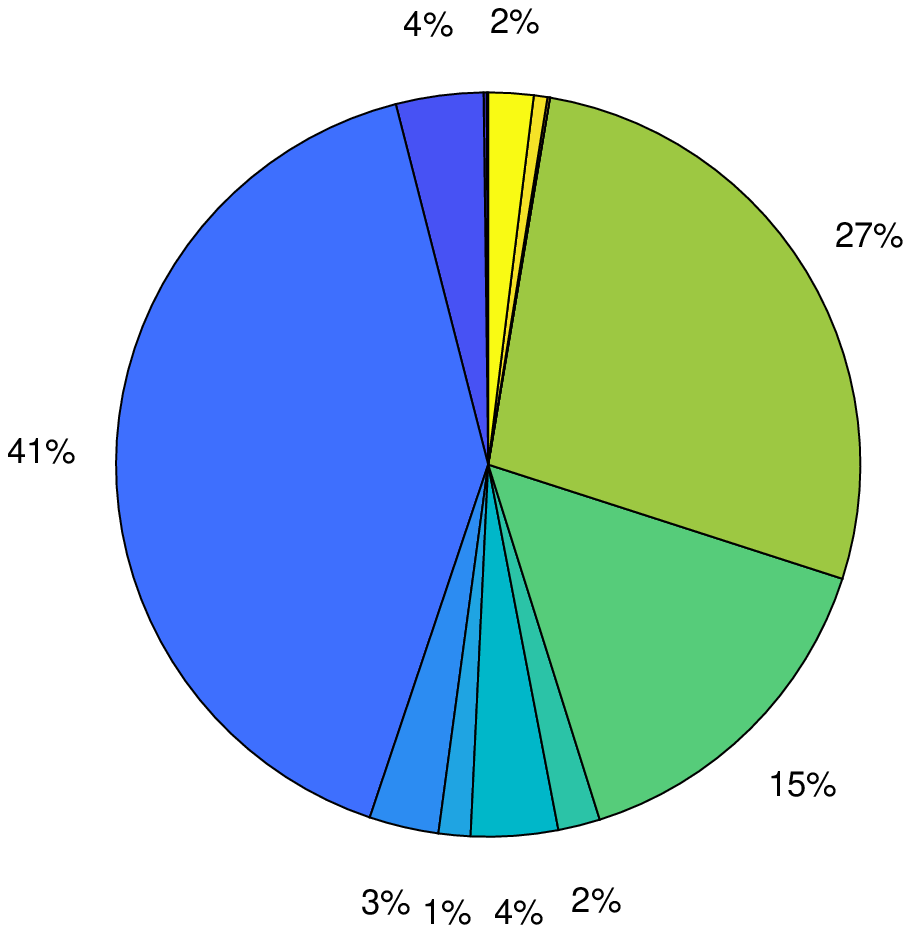}
\hspace{15mm}
\includegraphics[width=4.5cm]{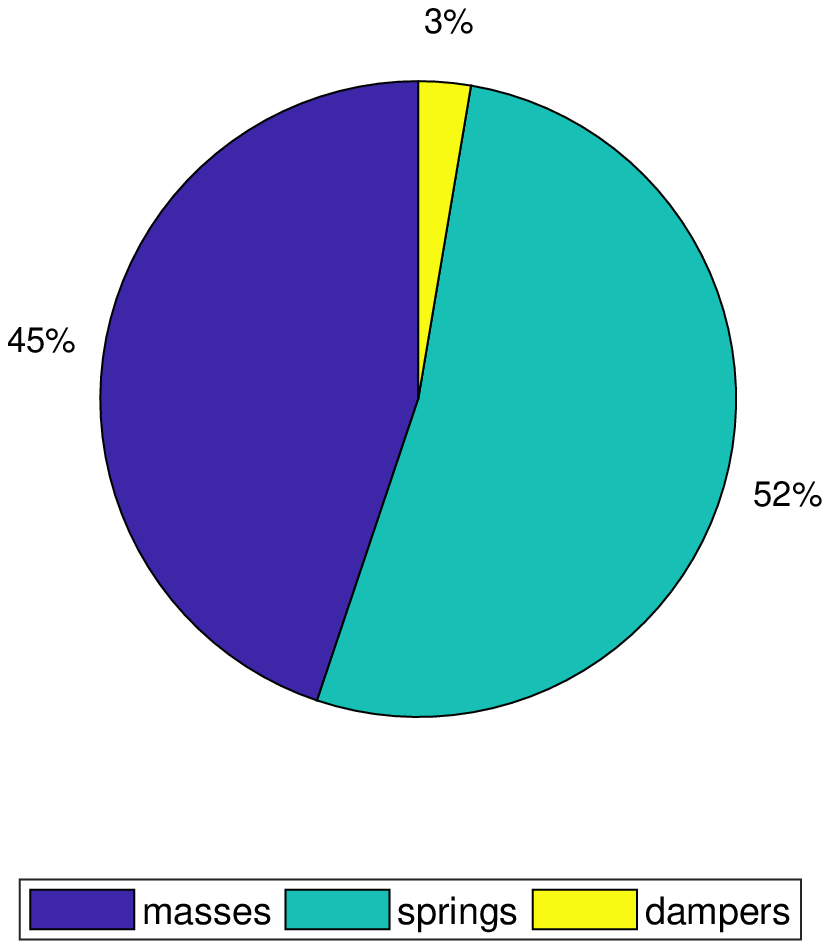}
\caption{Normalised sensitivity measures~(\ref{sensitivity-normalised2})
  for the random parameters (left) and merged normalised sensitivity measures
  for the three groups of physical parameters (right).}
\label{fig:pie}       
\end{figure}

The sensitivities indicate the impact of the random parameters
on the random QoI.
The normalised sensitivity measures~(\ref{sensitivity-normalised2})
are illustrated by Figure~\ref{fig:pie} (left).
We observe that three parameters are dominant, which are
the top mass, the spring above the top mass and the spring connecting
bottom mass and top mass, cf. Figure~\ref{fig:msd}.
In contrast, five parameters exhibit a contribution less than 1\%.
Thus these parameters can be remodelled by constants to reduce the
dimension of the random space.
Furthermore, Figure~\ref{fig:pie} (right) shows the sum of the
normalised sensitivity measures~(\ref{sensitivity-normalised2}) for
the three different types of parameters.
We recognise that the dampers feature a relatively low impact on the
random QoI.

\subsection{Comparison of MOR methods}
We compare the error between FOM and ROM for the computed
sensitivity measures~(\ref{sens-measure-hankel}).
The error indicator is defined as
\begin{equation} \label{error-rom}
  \max_{j=1,\ldots,q} \left| \hat{\eta}_j - \tilde{\eta}_j \right| .
\end{equation}
We do not consider the errors of the sensitivity measures~(\ref{sens-measure}),
because the computation work for the system norms of the
FOMs~(\ref{galerkin2}) is too high to be performed on our computer.

The following three variants of MOR methods are applied in MATLAB:
\begin{itemize}
\item[i)]
  Balanced truncation with direct solution of the Lyapunov equations using
  the function {\tt lyapchol} of the  control system toolbox.
  In the direct method, we use the matrices 
  $\hat{A}'=\hat{E}^{-1} \hat{A}$, $\hat{B}'=\hat{E}^{-1} \hat{B}$
  and  $\hat{E}'=I$ (identity) to simplify the Lyapunov equations,
  because the numerical solution becomes significantly faster.
  These matrices can be computed with negligible effort
  by a sparse $LU$-decomposition of~$\hat{E}$.
\item[ii)]
  Balanced truncation with iterative solution of the Lyapunov equations
  by the ADI method using the function {\tt lyapchol} of the sss toolbox,
  see~\cite{castagnotto}.
  The iteration performs until a default tolerance is achieved,
  which yields approximate factors of rank~301.
\item[iii)]
  One-sided Arnoldi method, see~\cite[p.~334]{antoulas},
  with a single real expansion point~$s_0=0$.
  The projection matrix~$T_{\rm r}$ is determined by the matrices
  $\hat{A},\hat{B},\hat{E}$ of~(\ref{galerkin}) and thus this approach is
  independent of the definition of outputs.
  It holds that $T_{\rm l}=T_{\rm r}$ used in~(\ref{reduced-matrices}).
\end{itemize}

In each technique, we calculate projection matrices $T_{\rm l},T_{\rm r}$
of rank~$r_{\max}=100$.
We extract the leading columns of the projection matrices to obtain ROMs of
dimensions $r=5,10,15,\ldots,100$.
Furthermore, the method~(i) yields all Hankel singular values
in the decomposition~(\ref{svd}).

\begin{figure} 
\centering
\includegraphics[width=6.5cm]{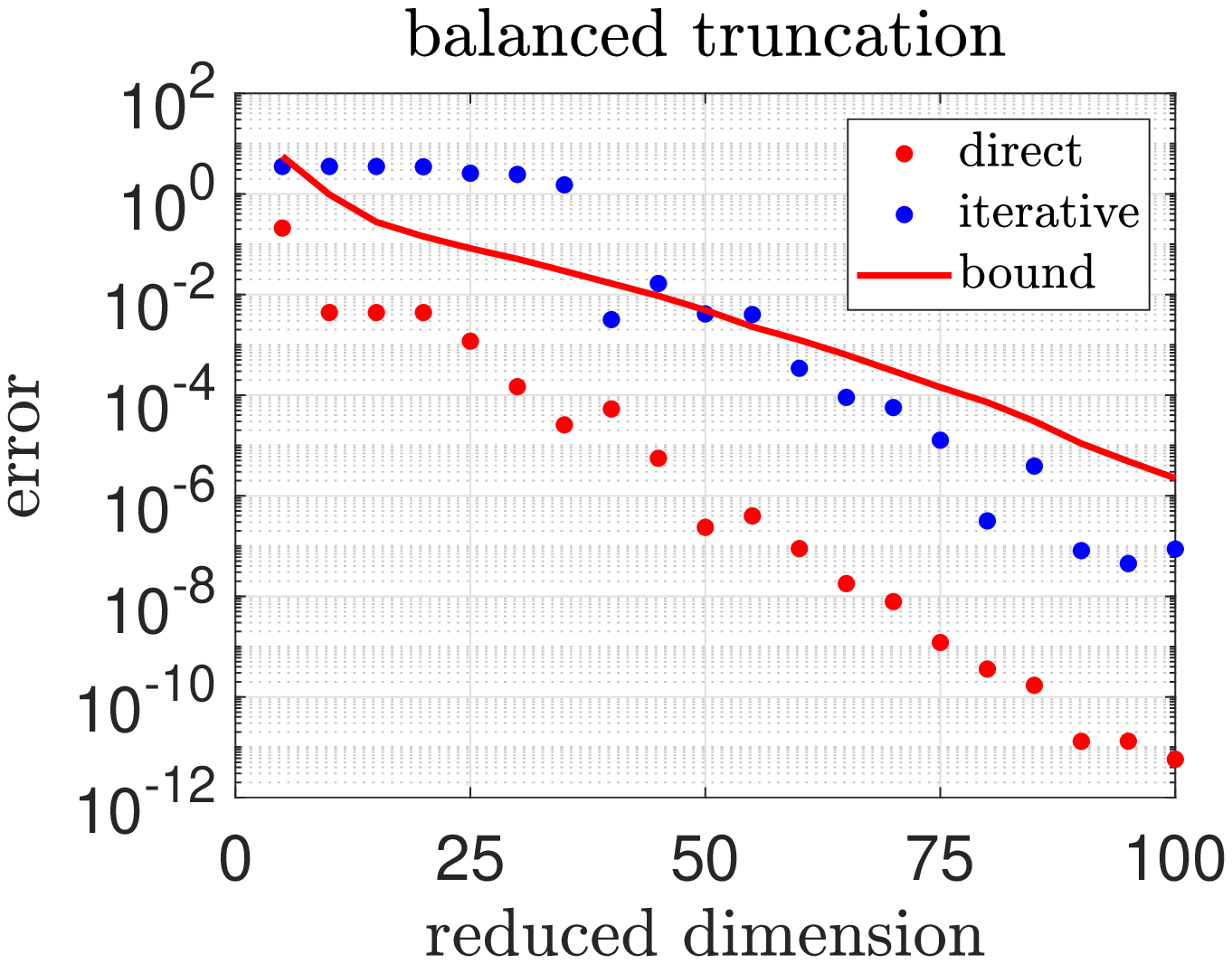}
\hspace{5mm}
\includegraphics[width=6.5cm]{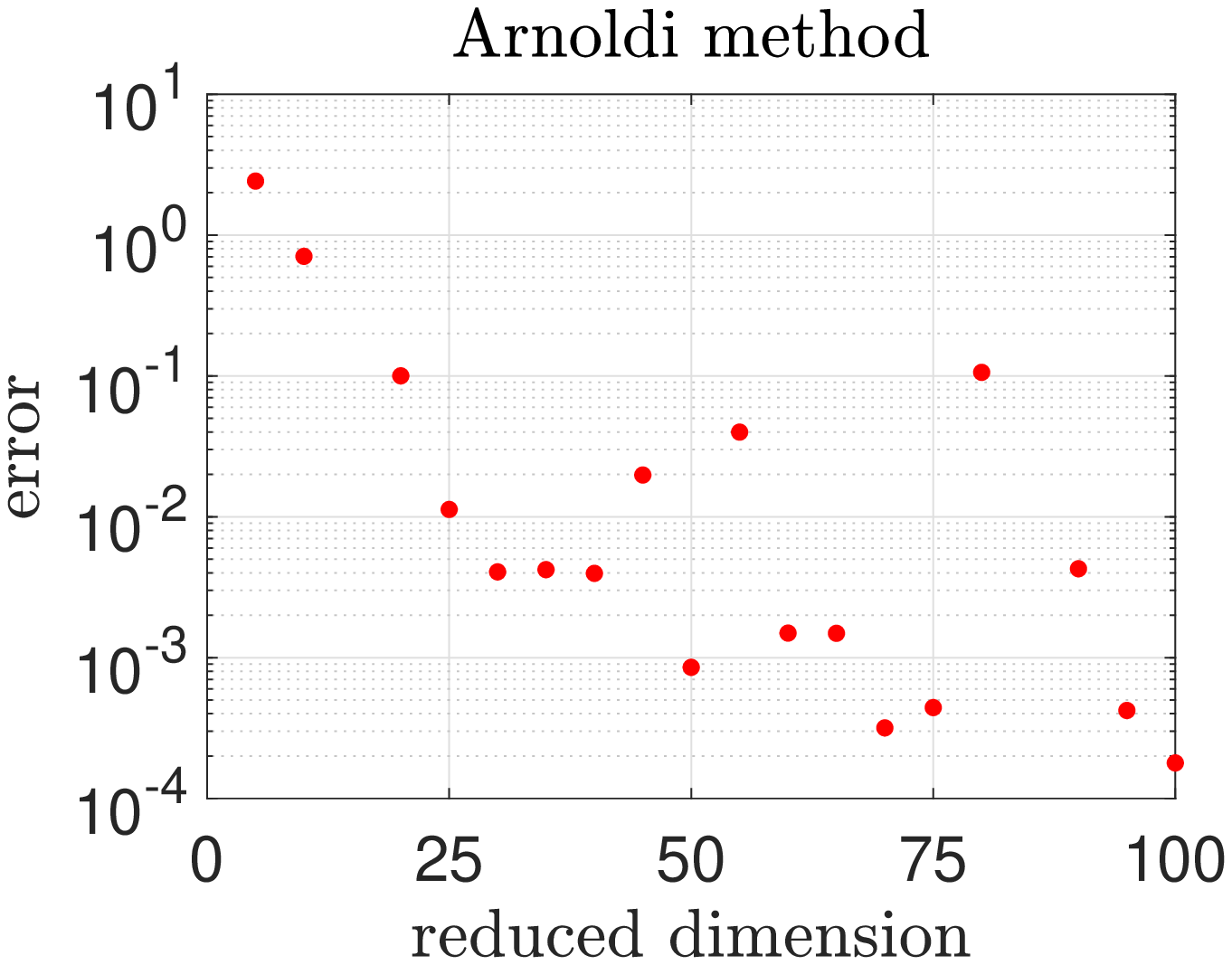}
\caption{Errors~(\ref{error-rom}) for different reduced dimensions
  in balanced truncation techniques (left) and Arnoldi method (right).}
\label{fig:errors}       
\end{figure}

Figure~\ref{fig:errors} (left) shows the maximum error~(\ref{error-rom})
for the techniques~(i) and~(ii).
Therein, the error bound~(\ref{error-bound}) obtained by
Theorem~\ref{thm:error} is also depicted.
In the direct method~(i), the true error is much below the
predicted error bound.
The approximate method~(ii) produces significantly higher errors
in comparison to the direct method,
which indicates a critical behaviour due to the large number of outputs
in the system~(\ref{galerkin}).
Nevertheless, the accuracy becomes better for increasing dimensions of
the ROMs.
Figure~\ref{fig:errors} (right) displays the maximum error~(\ref{error-rom})
for the technique~(iii). 
This error is larger for reduced dimensions $r \ge 50$ in comparison
to~(i) and~(ii).
Furthermore, the accuracy does not improve any more for dimensions
$r \ge 100$ in the Arnoldi method.
This is a well-known phenomenon due to the accumulation of
round-off errors in the orthogonalisation procedure.
Moreover, there is no prior or posterior error bound for the
Arnoldi method.

\begin{figure}
    \centering
    \includegraphics[width=6.5cm]{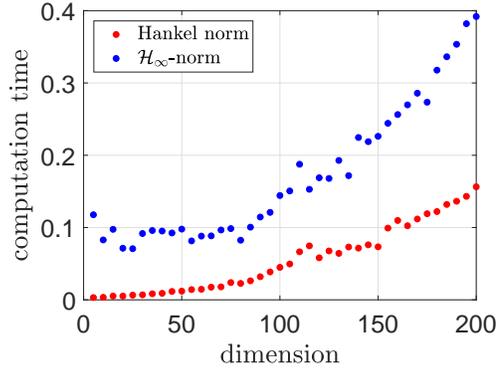}
    \caption{Computation times in seconds for a Hankel norm and a
      $\hinfty$-norm depending on the dimension of a linear dynamical system
      with dense matrices.}
    \label{fig:computation-times}
  \end{figure}

Table~\ref{tab:times} provides the computation times,
which represent elapsed real times in seconds, of the balanced
truncation techniques.
Only the $r_{\max}$ dominant singular values and their singular vectors
are computed in the SVD~(\ref{svd}) now.
The Arnoldi method requires a computation time of only 1.3 seconds for
the projection matrix of rank~$r_{\max}$.
Furthermore, the computation times for a single Hankel norm and
a single system norm of an ROM are illustrated by
Figure~\ref{fig:computation-times}.
In this case, we use the direct method~(i) to arrange ROMs up to
dimension $r=200$.
Yet the computational effort depends mainly on the reduced dimension and
not on the used MOR method.
The work has to be done for each random variable, i.e., $q$~times.
We recognise that this computation time is small in comparison to the
construction of the ROMs by balanced truncation techniques.

  \begin{table} 
    \caption{Computation times in seconds within
      balanced truncation MOR.\label{tab:times}}
    \begin{center}
      \begin{tabular}{lcc}
        & direct method & ADI method \\ \hline
        Lyapunov equation~(\ref{control-gramian}) & 467 & -- \\
        Lyapunov equation~(\ref{observ-gramian}) & 695 & -- \\
        both Lyapunov equations & 1162 & 26 \\
        singular value decomposition~(\ref{svd}) & 19 & 0.3
      \end{tabular}
    \end{center}
  \end{table}

The comparison demonstrates that there is a trade-off between accuracy and
computational effort in the MOR methods.
The balanced truncation using a direct solver is expensive,
whereas an a-priori error bound is guaranteed.
The direct method is still much cheaper than a computation of the
system norms for the FOMs~(\ref{galerkin2}).
  
\subsection{Transient simulation}
\begin{figure}
\centering
\includegraphics[width=6.5cm]{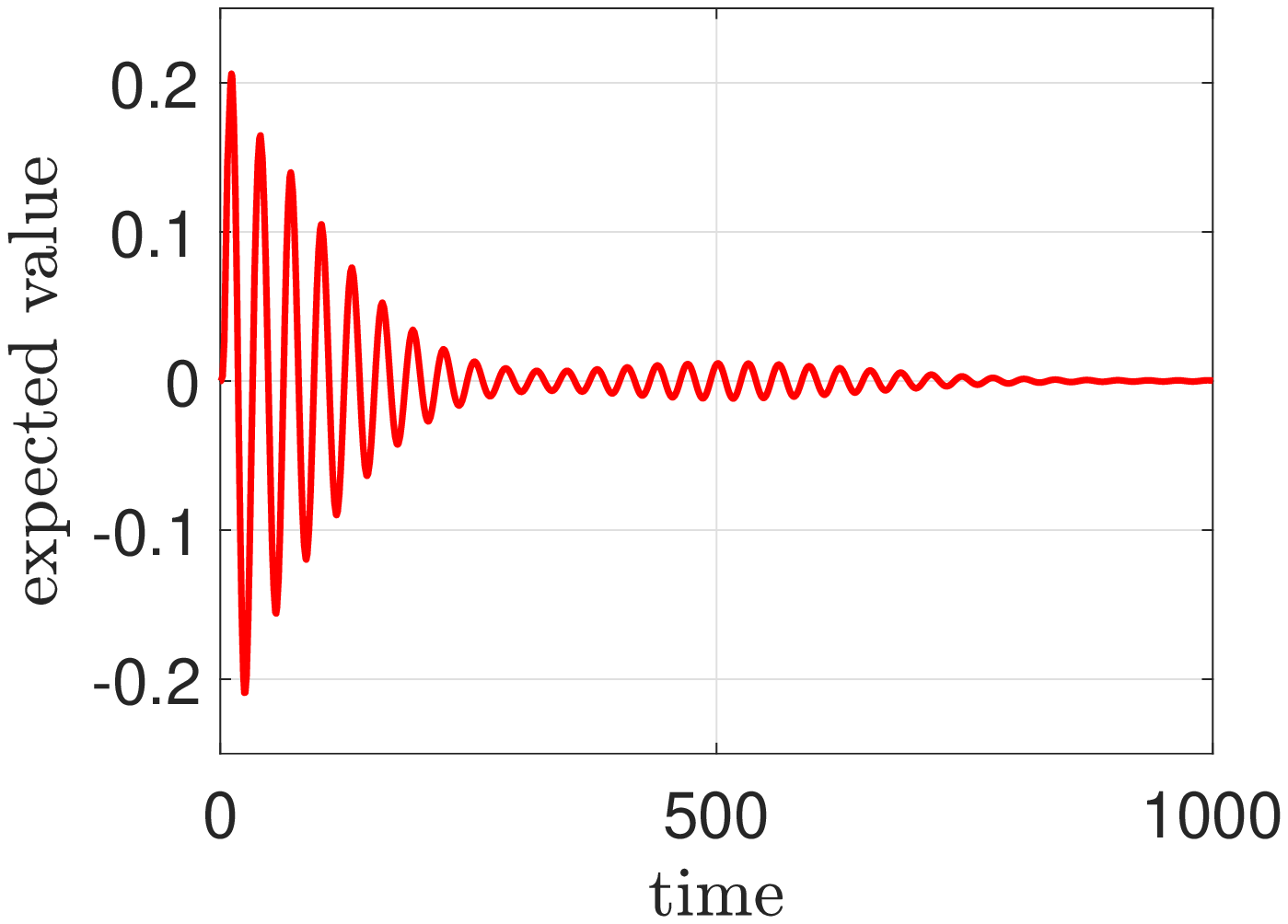}
\hspace{5mm}
\includegraphics[width=6.5cm]{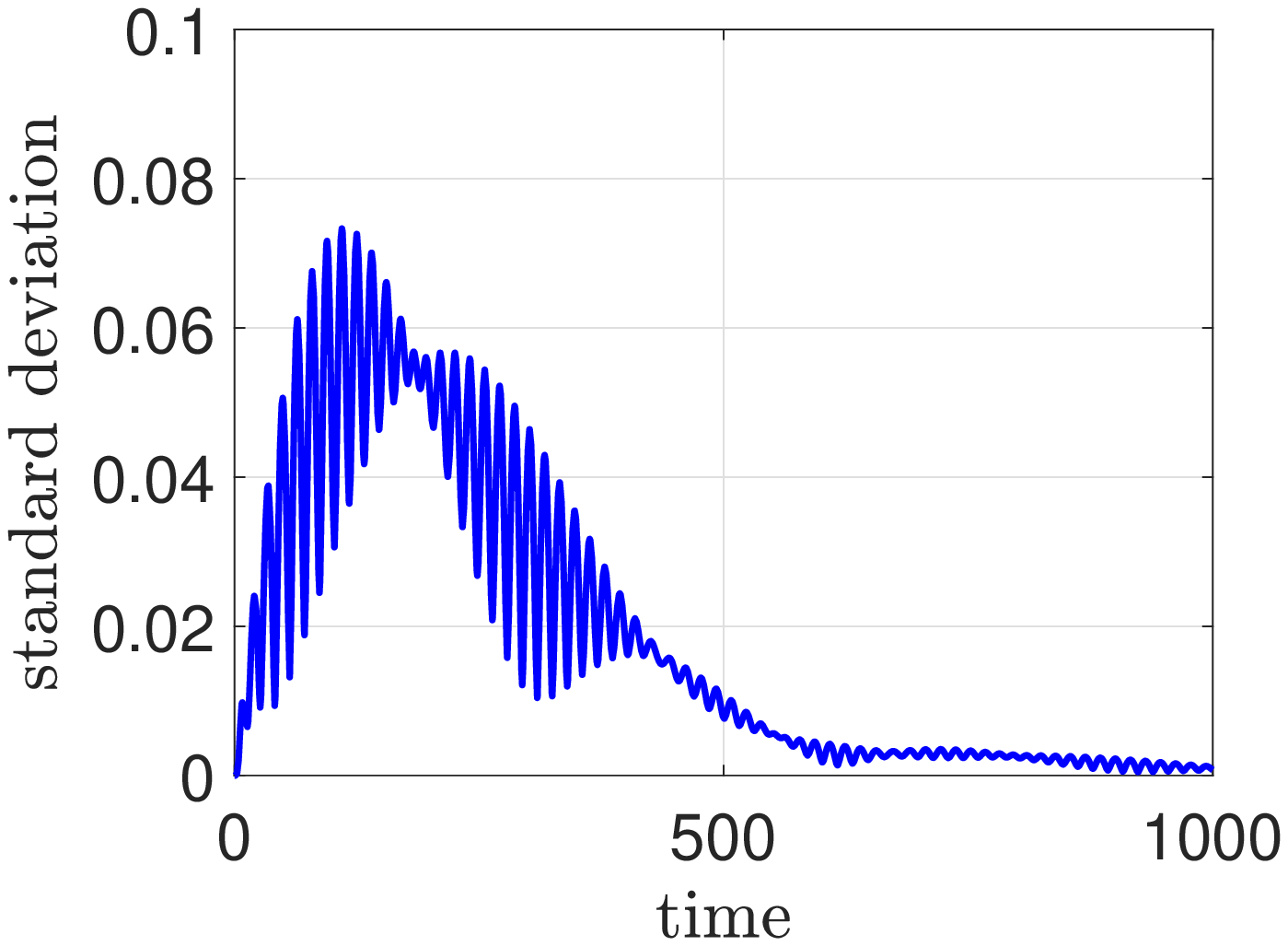}
\caption{Expected value (left) and standard deviation (right)
  of random output in transient simulation of mass-spring-damper system.}
\label{fig:transient}       
\end{figure}

\begin{figure}
\centering
\includegraphics[width=6.5cm]{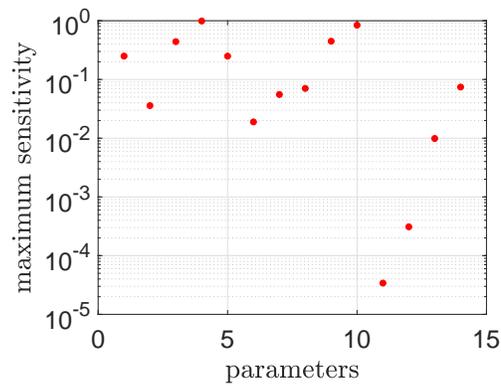}
\caption{Maximum in time of total effect sensitivity
  indices~(\ref{sens-appr}) 
  with respect to the random parameters
  (1-4: masses, 5-10: springs, 11-14: dampers).}
\label{fig:transient-sens}       
\end{figure}

Finally, we perform a transient simulation of the Galerkin
system~(\ref{galerkin}) for $t \in [0,1000]$ with initial values zero.
The chirp signal
$$ u(t) = \sin \left( \textstyle \frac{t^2}{10} \right) $$
is supplied,
which runs through a continuum of frequencies.
We use the trapezoidal rule in the time integration.
Figure~\ref{fig:transient} depicts the approximations of
the expected value and the standard deviation for the random QoI.
We compute the approximations~(\ref{sens-appr}) of the total effect
sensitivity indices~(\ref{total-sensitivity}) for each random variable.
Concerning the discrete time points, 
the maximum values are illustrated by Figure~\ref{fig:transient-sens}.
The relative positions of the sensitivity indices are similar to
the sensitivity measures~(\ref{sens-measure-hankel}),(\ref{sens-measure})
shown in Figure~\ref{fig:sensitivities}.

\clearpage


\end{document}